\documentclass[12pt]{amsart}
\usepackage{amsmath,amscd,amssymb,amsfonts,graphics,mathrsfs}
\usepackage{hyperref}
\newcommand{\hl}{\hyperlink}
\newcommand{\htt}{\hypertarget}
\newcommand{\h}{\hbox}
\newcommand{\q}{\quad}
\newcommand{\nin}{\noindent}
\newcommand{\bs}{\par\bigskip}

\newcommand{\sk}{\par\smallskip}
\newcommand{\bsn}{\par\bigskip\noindent}
\newcommand{\msn}{\par\medskip\noindent}

\newcommand{\ges}{\geqslant}

\newcommand{\1}{\hskip1pt}
\newcommand{\mcap}{\hbox{$\bigcap$}}
\newcommand{\mcup}{\hbox{$\bigcup$}}
\newcommand{\mopl}{\hbox{$\bigoplus$}}
\newcommand{\msum}{\hbox{$\sum$}}

\newcommand{\A}{{\mathscr A}}
\newcommand{\B}{{\mathscr B}}
\newcommand{\D}{{\mathscr D}}
\newcommand{\Fc}{{\mathscr F}}

\newcommand{\Hc}{{\mathscr H}}
\newcommand{\Lc}{{\mathscr L}}
\newcommand{\OO}{{\mathscr O}}

\newcommand{\ut}{\widetilde{u}}
\newcommand{\Xt}{\widetilde{X}}
\newcommand{\Lh}{\widehat{\mathscr L}}

\newcommand{\M}{{\mathcal M}}
\newcommand{\RR}{{\mathbf R}}
\newcommand{\R}{{\mathbb R}}
\newcommand{\PP}{{\mathbb P}}
\newcommand{\Q}{{\mathbb Q}}
\newcommand{\C}{{\mathbb C}}
\newcommand{\N}{{\mathbb N}}

\newcommand{\Z}{{\mathbb Z}}
\newcommand{\ub}{\overline{u}}

\newcommand{\xib}{\overline{\xi}}
\newcommand{\taub}{\overline{\tau}}
\newcommand{\zb}{\overline{z}}

\newcommand{\Gr}{{\rm Gr}}

\newcommand{\al}{\alpha}
\newcommand{\be}{\beta}

\newcommand{\la}{\lambda}

\newcommand{\Om}{\Omega}
\newcommand{\De}{\Delta}

\newcommand{\dd}{\partial}
\newcommand{\ddd}{{\rm d}}
\newcommand{\eq}{\,{=}\,}
\newcommand{\gess}{\,{\ges}\,}

\newcommand{\sgt}{\,{>}\,}
\newcommand{\slt}{\,{<}\,}
\newcommand{\nes}{\,{\ne}\,}

\newcommand{\pl}{\1{+}\1}

\newcommand{\bl}{\bigl}
\newcommand{\br}{\bigr}
\newcommand{\sst}{\,{\subset}\,}
\newcommand{\stm}{\,{\setminus}\,}

\newcommand{\ins}{\,{\in}\,}
\newcommand{\tos}{\,{\to}\,}
\newcommand{\defs}{\,{:=}\,}
\newcommand{\col}{\,{:}\,}
\newcommand{\ssc}{\,\raise.15ex\hbox{${\scriptstyle\circ}$}\,}
\newcommand{\ssb}{\raise.15ex\h{${\scriptscriptstyle\bullet}$}}
\newcommand{\into}{\hookrightarrow}
\newcommand{\onto}{\twoheadrightarrow}
\newcommand{\simto}{\,\,\rlap{\hskip1.3mm\raise1.4mm\hbox{$\sim$}}\hbox{$\longrightarrow$}\,\,}
\begin{document}
\h{}\bs\bs 
\centerline{\large Limits of Hodge structures with quasi-unipotent monodromies}
\bs
\centerline{Morihiko Saito}
\bs\bs
\vbox{\nin\narrower\smaller
{\bf Abstract.} We survey a theory of limits of polarizable variations of real Hodge structure in the quasi-unipotent monodromy case using the V-filtration of Kashiwara and Malgrange indexed by rational numbers, which does not necessarily seem familiar to many people.}
\bs\bs
\centerline{\bf Introduction}
\msn
We know very well about the limits of polarizable variations of real Hodge structure on a punctured disk in the {\it unipotent\1} monodromy case, see \cite{Sch} (and \cite{St1}, \cite{GN} for the geometric case). This implies certain information in the {\it quasi-unipotent\1} monodromy case using the {\it unipotent base change.} However, the precise formulation does not necessarily seem familiar to many people, since the $V$\!-filtration of Kashiwara \cite{Ka} and Malgrange \cite{Ma} {\it indexed by\1} $\Q$ plays an essential role. (Note that the latter filtration was indexed by $\Z$ in the original papers, see also \cite{sup}, \cite[(3.4)]{gm1}.) We explain some details about these. (This should have been written more than 30 years ago.)
\bs\bs
\centerline{\bf 1.~Abstract case}
\bsn
In this section we explain the case of abstract polarizable variations of real Hodge structure.
\msn
{\bf 1.1.~Deligne extensions.} Let $\bl((\Lc,F),L\br)$ be a polarizable variation of $\R$-Hodge structure of weight $w$ on a punctured open disk $S^*\defs S\stm\{0\}$ with $S\sst\C$ an open disk. Here $(\Lc,F)$ is a filtered locally free sheaf with an integrable connection $\nabla$, $L$ is an $\R$-local system with an isomorphism $\C{\otimes}_{\R}L\eq{\rm Ker}\,\nabla\sst\Lc$, and for any $c\ins S^*$, $\bl((\Lc_{(c)},F),L_c\br)$ is an $\R$-Hodge structure of weight $w$ with
$$\Lc_{(c)}\defs\C{\otimes}_{\OO_{S^*\!\!,c}}\Lc_c=\Lc_c/(t{-}c)\Lc_c.$$
We assume that it is {\it polarizable\1} by an pairing of the local system $L$, and the monodromy of $L$ is {\it quasi-unipotent.}
\sk
Let $I\sst\R$ be an {\it interval,} that is, a {\it connected\1} subset, satisfying the condition:
\htt{1.1.1}{}
$$\h{$I\,{\ni}\,\al\,\mapsto\,e^{2\pi i\al}\ins\{\la\ins\C\mid|\la|\eq1\}\,$ is bijective.}
\leqno(1.1.1)$$
Let $t$ be the coordinate of $S\sst\C$ (which is identified with the last inclusion). Let $j\col S^*\,{\into}\,S$ be the inclusion. Set $L_{\infty}\defs\Gamma(\widetilde{S}^*,\rho^*L)$ with $\rho\,{:}\;\widetilde{S}^*\tos S^*$ a universal covering. The following is due to \cite{De1} except the assertion on the $\Gr_F^p\Lc^I$ which follows from the {\it nilpotent orbit theorem\1} \cite{Sch}.
\par\htt{T1.1}{}\msn
{\bf Theorem~1.1}\,(Regularity). {\it There is uniquely a locally free extension $\Lc^I$ of $\Lc$ over $S$ such that the $\Gr_F^p\Lc^I$ are locally free $(\forall\,p\ins\Z$), $\Lc^I\sst j_*\Lc$ is stable by $t\dd_t\defs t\nabla_{\!\dd_t}$, and the eigenvalues of the residue ${\rm Res}^It\dd_t\ins{\rm End}_{\C}(\Lc^I_{(0)})$ are contained in $I$.}
\msn
{\it Proof.} There is the locally free extension $\Lc^I$ satisfying the condition of the connection by \cite{De1}, see also Remark~\hl{R1.1}{1.1} below. As for the filtration $F$, it is enough to show that each $F^p\Lc$ can be extended to a {\it coherent\1} subsheaf $\Fc^p\sst\Lc^I$. Indeed, we have the inclusion
\htt{1.1.2}{}
$$(\Lc^I\cap j_*F^p\Lc)/\Fc^p=\Gamma_{\{0\}}(\Lc^I/\Fc^p)\into\Lc^I/\Fc^p.
\leqno(1.1.2)$$
So the left-hand side is coherent, since $\Lc^I/\Fc^p$ is. The intersection $\Lc^I\cap j_*F^p\Lc$ is thus coherent, and we can put
\htt{1.1.3}{}
$$\Fc^p\defs\Lc^I\cap j_*F^p\Lc.
\leqno(1.1.3)$$
\sk
To get a coherent extension $\Fc^p$, we may assume that $I\eq[0,1)$ (taking the intersection of $\Fc^p$ with $\Lc^I$ if $I\nes[0,1)$). Let $\rho'\;{:}\,S'\tos S$ be a ramified cyclic covering of degree~$m$ such that $L'\defs\rho'{}^*L$ has unipotent monodromy, where $S'$ is an open disk. Let $t'$ be a coordinate of $S'$ such that $\rho'\1^*t\eq t'\1^m$. Let $T,T'$ be the monodromies of $L,L'$ with Jordan decompositions $T\eq T_sT_u$, $T'\eq T'_sT'_u$ respectively. Put $N\defs\log T_u$, $N'\defs\log T'_u$. Then $T'$ is identified with $T\1^m$, and $N'$ with $m\1N$ via the isomorphism $L_{\infty}\eq L'_{\infty}$.
\sk
Let $u_1,\dots,u_r$ be a basis of $L_{\C,\infty}$ such that $T_su_j\eq\la_ju_j$ ($\la_j\ins\C$), where $L_{\C}\defs\C{\otimes}_{\R}L$. Let $\ut_1,\dots,\ut_r$ and $\ut'_1,\dots,\ut'_r$ be the corresponding free generators of the Deligne extensions $\Lc^I$ on $S$ and $\Lc'{}^{[0,1)}$ on $S'$ respectively, see Remark~\hl{R1.1}{1.1} below. Take $m_j\ins\Z$ such that $\al_j\defs\tfrac{m_j}{m}\ins[0,1)$ and $e^{-2\pi i\al_j}\eq\la_j$. By construction (see (\hl{1.1.6}{1.1.6}) below), we have
\htt{1.1.4}{}
$$\rho'{}^*\ut_j\eq t'{}^{m_j}\ut'_j\q(j\ins[1,r]).
\leqno(1.1.4)$$
Moreover the $t'{}^i\ut'_j$ ($i\ins[0,m{-}1],\,j\ins[1,r]$) are free generators over $\C\{t\}\,({\subset}\,\C\{t'\})$ of $\Lc^{\prime[0,1)}_0$ (which is identified with $(\rho'_*\Lc'{}^{[0,1)})_0$), and $m_j\ins[0,m{-}1]$ ($j\in[1,r]$). We then see that there is a surjective morphism of $\C\{t\}$-modules
\htt{1.1.5}{}
$$\Lc^{\prime[0,1)}_0\onto\Lc^{[0,1)}_0,
\leqno(1.1.5)$$
splitting the inclusion $\rho'{}^*\col\Lc^{[0,1)}_0\into\Lc^{\prime[0,1)}_0$. So the assertion is reduced to the unipotent monodromy case, and follows from the {\it nilpotent orbit theorem\1} \cite{Sch}. This finishes the proof of Theorem~\hl{T1.1}{1.1}.
\par\htt{R1.1}{}\msn
{\bf Remark 1.1} (Deligne extension \cite{De1}). For a $\C$-local system $L_{\C}$ on $S^*$, set $\Lc\defs\OO_{S^*}{\otimes}_{\C}L_{\C}$. This has the natural connection $\nabla$. Let $I\sst\R$ be any subset satisfying (\hl{1.1.1}{1.1.1}). Then there is uniquely a locally free extension $\Lc^I$ of $\Lc$ over $S$ in $j_*\Lc$ such that $\Lc^I$ is stable by $t\dd_t\defs t\nabla_{\!\dd_t}$, and the eigenvalues of the residue ${\rm Res}^It\dd_t\ins{\rm End}_{\C}(\Lc^I_{(0)})$ are contained in $I$. Indeed, for a basis $u_1,\dots,u_r$ of $L_{\C,\infty}$ with $T_su_j\eq\la_ju_j$ ($\la_j\ins\C^*,\,j\ins[1,r])$, set
\htt{1.1.6}{}
$$\ut_j\defs t^{\al_j}e^{-(2\pi i)^{-1}(\log t)N}u_j\in j_*\Lc\q(\al\in I,\,e^{-2\pi i\al_j}\eq\la_j).
\leqno(1.1.6)$$
Then $\ut_1,\dots,\ut_r$ are {\it free generators of the Deligne extension\1} $\Lc^I\sst j_*\Lc$, and there is an isomorphism of $\C$-vector spaces
\htt{1.1.7}{}
$$L_{\C,\infty}=\Lc^I_{(0)}\,(=\Lc^I_0/t\Lc^I_0),
\leqno(1.1.7)$$
such that the action of the monodromy $T$ of the local system on the left-hand side corresponds to $e^{-2\pi i\,{\rm Res}^It\dd_t}$ on the right-hand side
\msn
{\bf 1.2.~Shifted $V$-filtration.} In the notation of 1.1, set for $\al\ins\Q$
$$\aligned\Lc^{\ges\al}\defs\Lc^I\q\h{if}\,\,\,I\eq[\al,\al{+}1),\\ \Lc^{>\al}\defs\Lc^I\q\h{if}\,\,\, I\eq(\al,\al{+}1].\endaligned$$
By the construction of $\Lc^I$ (see (\hl{1.1.6}{1.1.6})), we have the inclusions
\htt{1.2.1}{}
$$\Lc^{\ges\al}\,{\supset}\,\Lc^{>\al}\,{\supset}\,\Lc^{\ges\be}\q\h{if}\,\,\,\al\slt\be.
\leqno(1.2.1)$$
Put
$$\Lh\defs\mcup_{\al\in\Q}\,\Lc^{\ges\al}\eq\mcup_{\al\in\Q}\,\Lc^{>\al}.$$
This is a {\it meromorphic regular holonomic\1} $\D_S$-module (so the action of $t$ is bijective). Set
$$V'{}^{\al}\Lh\defs\Lc^{\ges\al},\q V'{}^{>\al}\Lh\defs\Lc^{>\al},\q\Gr_{V'}^{\al}\Lh\defs V'{}^{\al}\Lh/V'{}^{>\al}\Lh.$$
This coincides, up to the shift of index by 1, with the $V$\!-filtration of Kashiwara \cite{Ka} and Malgrange \cite{Ma} {\it indexed by\1} $\Q$. Note that it was indexed by $\Z$ in the original papers where no Hodge filtration appears, see also \cite{sup}, \cite[(3.4)]{gm1}, and \cite{Va}. (Perhaps $V$ may come from the latter.)
\sk
Using ${\rm Ker}\,(t\dd_t{-}\al)^r\sst\Lh$ for $\al\ins\Q$ with $r\eq{\rm rank}\,L$, the filtration $V'$ on $\Lh$ splits canonically (fixing $t$), and we get the isomorphism
\htt{1.2.2}{}
$$\Lc^I_{(0)}=\mopl_{\al\in I}\,\Gr_{V'}^{\al}\Lh_0,
\leqno(1.2.2)$$
corresponding to the decomposition by monodromy eigenvalues of the left-hand side of (\hl{1.1.7}{1.1.7}). We have the filtration $F$ on the right-hand side of (\hl{1.2.2}{1.2.2}) by
\htt{1.2.3}{}
$$F^p\Gr_{V'}^{\al}\Lh_0:=F^pV'{}^{\al}\Lh_0/F^pV'{}^{>\al}\Lh_0\into\Gr_{V'}^{\al}\Lh_0,
\leqno(1.2.3)$$
which is {\it independent\1} of $I$ containing $\al$, where $\Lh\sst j_*\Lc$ has the filtration $F$ defined by
$$F^p\Lh\defs\Lh\cap j_*F^p\Lc\q(p\ins\Z).$$
(This is {\it different\1} from the Hodge filtration of mixed Hodge module.)
\sk
Let $W$ be the {\it monodromy filtration\1} on $L_{\infty}$ shifted by $w$. This is uniquely characterized by the two conditions:
\htt{1.2.4}{}
$$\aligned&\q\q N(W_kL_{\infty})\sst W_{k-2}L_{\infty}\q(\forall\,k\ins\Z),\\&N^k:\Gr^W_{w+k}L_{\infty}\simto\Gr^W_{w-k}L_{\infty}\q(\forall\,k\ins\N).\endaligned
\leqno(1.2.4)$$
This induces the filtration $W$ on $\Lc^I_{(0)}$ via the isomorphism (\hl{1.1.7}{1.1.7}).
\par\htt{T1.2}{}\msn
{\bf Theorem~1.2}\,(Limit mixed Hodge structure). {\it The regular extension in Theorem~{\rm \hl{T1.1}{1.1}} gives a limit mixed $\R$-Hodge structure $\bl((\Lc^I_{(0)};F,W),(L_{\infty},W)\br)$ compatible with the one constructed by Schmid in the unipotent monodromy case via a unipotent base change as in the proof of Theorem~{\rm \hl{T1.1}{1.1}}, where $L_{\infty}$ does not change by the unipotent base change.}
\msn
{\it Proof.} The limit Hodge filtration on $\Gr_{V'}^{\al}\Lh$ is independent of $I$, see (\hl{1.2.3}{1.2.3}). The assertion is then reduced to the case $I\eq[0,1)$ in view of (\hl{1.1.6}{1.1.6}), and follows from Schmid's SL$_2$-orbit theorem \cite{Sch} in the unipotent monodromy case using (\hl{1.1.4}{1.1.4}). Note that the weight filtration is also given by the shifted monodromy filtration in Schmid's theorem. This finishes the proof of Theorem~\hl{T1.2}{1.2}.
\par\htt{R1.2a}{}\msn
{\bf Remark~1.2a.} Theorem~\hl{T1.2}{1.2} does not hold if $I$ is {\it not connected,} that is, if the Deligne extension is {\it not\1} a member of the shifted $V$-filtration. Consider for instance the case $L_{\C,\infty}$ is generated by $u$ and $u'\defs\ub$ with $Tu=\la u$ ($\la\ins\C$), and the Hodge filtration $F^1\Lc$ is generated by
\htt{1.2.5}{}
$$v\defs t^{\al}u\pl t^{\al'}u'\q\h{with}\q\la\eq e^{-2\pi i\al},\,\al\ins\bl(0,\tfrac{1}{2}\br),\,\al{+}\al'\eq1,
\leqno(1.2.5)$$
where the weight of the variation of Hodge structure is 1. This gives an example such that the ``naive" limit Hodge filtration on $\Lc^I_{(0)}$ for $I\eq[0,1)$ is not compatible with the decomposition by monodromy eigenvalues, since $\al,\al'\ins[0,1)$, see \cite{sup}.
\sk
Moreover, if we take $I$ such that $\al\1{+}\1i\ins I$ and $\al'{+}\1i'\ins I$ for $i,i'\ins\Z$ with $i\slt i'$ (where $I$ cannot be \h{{\it connected}\1)}, then we see that the ``leading term" of the limit in $\Lc^I_{(0)}$ is given by $[t^{\al'+i'}u']$ rather than $[t^{\al+i}u]$, considering $[t^{i'}v]\ins\Lc^I_{(0)}$, where $t^iv\,{\notin}\,\Lc^I$.
\par\htt{R1.2b}{}\msn
{\bf Remark~1.2b.} The eigenvalues of the monodromy have absolute value 1 even in the case of polarizable variations of {\it complex\1} Hodge structure, see \cite{SabSc}. (It is easy to verify this in the rank 1 case using the polarizability.) So the $V$\!-filtration is indexed by $\R$.
\bs\bs
\vbox{\centerline{\bf 2. Geometric case}
\bsn
In this section we explain geometric proofs of Theorems~\hl{T1.1}{1.1} and \hl{T1.2}{1.2} in the geometric case.}
\msn
{\bf 2.1.~Steenbrink's method in the reduced case \cite{St1}.} Let $f\col X\tos S$ be a proper morphism from a connected K\"ahler manifold onto an open disk $S$. We may assume that $f$ induces a smooth morphism over $S^*\defs S\stm\{0\}$ and that the singular fiber $Y\defs f^{-1}(0)$ is a divisor with simple normal crossings using a desingularization \cite{Wlo}, since we are interested in the limit mixed Hodge structures. In this subsection we assume
\htt{2.1.1}{}
$$\h{The singular fiber $Y$ is {\it reduced.}}
\leqno(2.1.1)$$
This is satisfied in the {\it semistable reduction\1} case \cite{KKMS}. (It does not seem very clear whether the argument in the algebraic case applies entirely in the same way to the analytic case, since some people seem to say that there is some problem.)
In this case Theorem~\hl{T1.1}{1.1} and \hl{T1.2}{1.2} can be proved by an essentially same argument as in \cite{St1} using \cite[5.12]{StZ} for the double complex construction with $\R$-coefficients and \cite[4.2.2]{mhp} for the {\it stability of the monodromy filtration\1} under the direct image by $f$, see also \cite{GN}, \cite{mon}.
\sk
Note that the double complex construction in \cite[(4.14)]{St1} gives an explicit description of the {\it monodromy filtration\1} on the nearby cycle complex $\psi_f\C_X[n]$, where $n\eq\dim Y$.
\par\htt{R2.1a}{}\msn
{\bf Remark~2.1a.}  It is quite surprising that a presented {\it correction argument\1} contains a serious {\it confusion\1} between ``local" in SGA7, XIV, 4.18 and ``global" in ASENS 19, p.\,127 (see also \cite[2.4.2]{Il}); for instance, consider the case where the irreducible components of $Y$ are {\it rational\1} surfaces and some of their intersections is a {\it non-rational\1} curve; more concretely, the blow-up at $(0,0)$ of the family $\{f_3w{+}f_4\eq tw^4\}\sst\PP^3{\times}\C$, where $f_k\ins\C[x,y,z]$, $\deg f_k\eq k$ $(k\eq3,4$), and $V_{f_3}\defs\{f_3\eq0\}$ is smooth and intersects transversally $V_{f_4}$ at 12 smooth points; it is also possible to take the blow-up at 16 singular points of $\{g\1g'\eq t\1h\}\sst\PP^3{\times}\C$, where $V_g,V_{g'}$ are nonsingular conics intersecting transversally and $V_g\1{\cap}\1V_{g'}$ intersects transversally $V_h$ at 16 smooth points (since $\deg h\eq4$).
\par\htt{R2.1b}{}\msn
{\bf Remark~2.1b.} We have to study the {\it weight spectral sequence\1} to prove the {\it hard Lefschetz\1} property for the action of $N$ on the $E_2$-term (that is, the {\it stability\1} of the monodromy filtration by the direct image) and also the {\it strictness\1} of the Hodge filtration $F$ on
\htt{2.1.2}{}
$$\RR\Gamma\bl(Y,\Om^{\ssb}_{X/S}(\log Y){\otimes}_{\OO_X}\OO_Y\br),
\leqno(2.1.2)$$
induced by the so-called ``stupid" filtration $\sigma_{\ges\ssb}$. The strictness of $F$ is needed to apply the {\it semicontinuity\1} argument for the proof of the {\it freeness\1} of the higher direct image sheaves
$$R^q\!f_*\Om^p_{X/S}(\log Y)\q(p,q\ins\N).$$
Here the assumption (\hl{2.1.1}{2.1.1}) {\it cannot\1} be replaced with the condition that the monodromy is {\it unipotent.} Indeed, to apply the semi-continuity argument, the filtered complex in (\hl{2.1.2}{2.1.2}) cannot be replaced by
$$\RR\Gamma\bl(Y,\Om^{\ssb}_{X/S}(\log Y){\otimes}_{\OO_X}\OO_{Y_{\rm red}}\br),$$
even though the two complexes are quasi-isomorphic by the assumption on the monodromy forgetting the filtration $F$, see also \cite[Remark after Cor.\,5.10]{St1}.
\msn
{\bf 2.2.~Steenbrink's method in the general case \cite{St2}.} Let $f\col X\tos S$ be as in 2.1 above, where (\hl{2.1.1}{2.1.1}) is not assumed. In this case we can take the {\it normalization of the base change\1} by a ramified cyclic covering $\rho'\col S'\tos S$ of degree $m$ as in the proof of Theorem~\hl{T1.1}{1.1} with $m$ the LCM of the multiplicities of $Y$ along irreducible components, see \cite{St2}. We then get V-manifolds (having quotient singularities), and the argument in 2.1 can be generalized.
\par\htt{R2.2}{}\msn
{\bf Remark~2.2.} This argument is used in \cite[(2.4)]{gm1} by considering
$$\A^{\ssb}\defs{\rm Ker}\bl(\ddd f{\wedge}\col\Om^{\ssb}_X\tos\Om^{\ssb+1}_X\br).$$
Via the morphism defined by $\ddd f\!/\!f{\wedge}$, this complex is filtered isomorphic to
$$\Om_{X/S}^{\ssb}(\log Y)(-Y_{\rm red})[-1].$$
\msn
{\bf 2.3.~Hodge module method \cite{mhp}, \cite{mhm}.} Let $f\col X\tos S$ be as in 2.1 with (\hl{2.1.1}{2.1.1}) not assumed. Let $i_f\col X\into X{\times}\C$ be the graph embedding. Set
$$(\B_f,F)\defs(i_f)_*^{\D}(\OO_X,F)\,\bl(=(\OO_X[\dd_t]\delta(t{-}f),F)\br).$$
Here $(i_f)_*^{\D}$ denotes the direct image for left $\D$-modules (where $F$ is shifted by the relative dimension 1), and $\Gr^F_p\OO_X\eq0$ ($p\nes 0$). It has the $V$\!-filtration of Kashiwara \cite{Ka} and Malgrange \cite{Ma} {\it indexed by\1} $\Q$ (where $\dd_tt{-}\al$ is nilpotent on $\Gr_V^{\al}$). The nearby and vanishing cycle functors are defined by
$$\aligned\psi_f(\OO_X,F)&\defs\mopl_{\al\in(0,1]}\,\Gr_V^{\al}(\B_f,F),\\ \varphi_{f,1}(\OO_X,F)&\defs\Gr_V^0(\B_f,F[-1]),\endaligned$$
where $\varphi_{f,1}$ denotes the unipotent monodromy part of the vanishing cycle functor $\varphi_f$. The non-unipotent monodromy part is defined by $\varphi_{f,\ne1}\eq\psi_{f,\ne1}$ with $\psi_f\eq\psi_{f,1}\oplus\psi_{f,\ne1}$. The weight filtration $W$ on $\psi_f\OO_X$ and $\varphi_{f,1}\OO_X$ is given by the monodromy filtration for the action of $N\defs\dd_tt{-}\al$, which is characterized as in (\hl{1.2.4}{1.2.4}) with $w\eq n$ and $n{+}1$ respectively. (Recall that $n\eq\dim Y$.).
\sk
Set $P_N\Gr^W_{w+k}\defs{\rm Ker}\,N^{k+1}\sst\Gr^W_{w+k}$ ($k\gess 0$) with $w\eq n$ or $n{+}1$. This is called the $N$-{\it primitive part.} It is relatively easy to see that the {\it unipotent\1} monodromy part is given by
\htt{2.3.1}{}
$$P_N\Gr^W_{n+k}\psi_{f,1}(\OO_X,F)=\mopl_{|I|=k+1}\,(i_{Y_I})_*^{\D}(\OO_{Y_I},F[-k]).
\leqno(2.3.1)$$
Here $Y_I\defs\mcap_{i\in I}\,Y_i$ with $Y_i$ the irreducible components of $Y$, and $i_{Y_I}\col Y_I\into X$ is the inclusion. (This is compatible with the calculation in \cite{St1}.) The description of the {\it non-unipotent\1} monodromy part is more complicated. Set $\psi_{f,\la}\defs{\rm Ker}(T_s{-}\la)\sst\psi_f$ ($\la\in\C$), and
$$Y_I^{(\la)}\defs Y_I\stm\mcup_{i\notin I_{\la}}\,Y_i\q\h{with}\q I_{\la}\defs\{i\,|\,\la^{m_i}\eq1\},$$
where $Y\eq\msum_i\,m_iY_i$. Then $P_N\Gr^W_{n+k}\psi_{f,\la}(\OO_X,F)$ corresponds via the de~Rham functor to the direct sum of the 0-extension of a unitary local system of rank 1 on $Y_I^{(\la)}$ for $I\sst I_{\la}$ with $|I|\eq k{+}1$. Here the eigenvalues of the local monodromies around the boundary divisors are not 1, and we get the {\it intermediate direct images\1} in the sense of \cite{BBD}. Using the local combinatorial description of regular holonomic $\D$-modules of normal crossing type (see for instance \cite[Rem.\,1.3b]{rh}), it can be described locally by a combinatorial formula, see \cite[Thm.\,3.3 and 3.4]{mhm}. We have the induced polarization using \cite{dual}. (To show that applying the nearby and vanishing cycle functors to $\bl((\OO_X,F),\RR_X[n{+}1]\br)$, we get mixed Hodge modules (in a weak sense) whose weight filtrations are shifted monodromy filtrations, it is enough to consider {\it locally\1} a ramified covering such that the pullback of $Y$ has {\it equal multiplicities\1} along irreducible components, where the calculation is essentially the same as in the reduced case, see \cite[5.4.3]{mhm}.)
\sk
We have the {\it weight spectral sequence\1} in the category of mixed $\R$-Hodge structures (which is identified with the category of mixed $\R$-Hodge modules on a point), and its filtered $\C$-vector space part is given by
\htt{2.3.2}{}
$$E_1^{-k,j+k}\eq H^j(a_X)_*^{\D}\Gr^W_k\psi_f(\OO_X,F)\Longrightarrow H^j(a_X)_*^{\D}\psi_f(\OO_X,F),
\leqno(2.3.2)$$
where $a_X\col X\tos pt$ denotes the structure morphism (similarly with $\psi_f$ replaced by $\varphi_{f,1}$). For the real coefficient part, we use Verdier's spectral object, see \cite[5.2.18]{mhp}. So we get the {\it strictness\1} of the Hodge filtration $F$ on $(a_X)_*^{\D}\psi_f(\OO_X,F)$ (using the abelian category of graded $\C$-vector spaces containing the exact category of finite-dimensional filtered $\C$-vector spaces as a full subcategory, see \cite[1.3.2]{mhp}) and also the {\it stability\1} of the monodromy filtration under the direct image by using \cite[4.2.2]{mhp} (see also \cite{mon}). This strictness implies the {\it bistrictness\1} of the direct image $(f{\times}id)_*^{\D}(\B_f;F,V)$ and the {\it stability\1} of the $V$-filtration by the direct image, see \cite[Prop.\,3.3.17]{mhp}. We then get the commutativity
\htt{2.3.3}{}
$$H^j(a_X)_*^{\D}\psi_f(\OO_X,F)=\psi_t\Hc^j\!f_*^{\D}(\OO_X,F)
\leqno(2.3.3)$$
(similarly with $\psi_f$ replaced by $\varphi_{f,1}$) and also the {\it freeness\1} of $F_pV^{>0}\Hc^j\!f_*^{\D}\OO_X$ (shrinking $S$) under the assumption that there is a cohomology class $\eta\ins H^2(Y,\R)$ whose restriction to any irreducible component of $Y$ is a {\it K\"ahler class.} We thus get a geometric proof of Schmid's theorem.
\par\htt{R2.3a}{}\msn
{\bf Remark~2.3a.} A {\it hybrid\1} argument is presented in \cite{Ch}, where one uses the {\it induced $\D_X$-module\1} $\M\defs{\rm DR}_X^{-1}\bl(\Om_{X/S}^{\ssb+n}(\log Y)|_Y\br)$ associated with the restriction to $Y$ (as $\OO$-module) of the relative logarithmic de~Rham complex in \cite{St1} instead of the nearby cycle $\D$-module (more precisely, $(V^0\!/V^1)\B_f(*Y)$), see \cite[2.2]{mhp}, \cite{ind} for the functor ${\rm DR}^{-1}_X$. These are locally isomorphic using the local combinatorial description of regular holonomic $\D$-modules of normal crossing type, see \cite[Rem.\,1.3b]{rh}, \cite[Thm.\,3.3--4]{mhm}. This is combined with the {\it hermitian pairings\1} defined by using the {\it residues of complex Mellin transforms\1} in \cite{Sab1}, \cite{Sab2} in order to get a {\it complex\1} mixed Hodge module {\it without\1} real structure. It may be interesting to compare the argument on the induced polarization in \cite{Ch} with the proof of \cite[Prop.\,3.8.1]{Sab2} where the residues of Mellin transforms are also used (although it contains a problem of {\it half Tate twists,} see \cite{twi}). The $V$\!-filtration of Kashiwara and Malgrange is {\it not\1} used in \cite{Ch}, although it appears implicitly in the calculation of the nearby cycle $\D$-module in the {\it non-reduced\1} case, where the $V$\!-filtration is essentially identified with the {\it multiplier ideals,} see \cite{BS}. (This is not surprising by the definition of the latter using the local integrability of $|g|^2/|f|^{2\al}$ for $g\in\OO_X$.)
\sk
One problem is that the description via the ideal of the $\D$-module $\M$ depends heavily on the choice of {\it local coordinates and their order\1} (although the formula itself is independent apparently), since the same applies to the corresponding {\it local generator\1} which is annihilated by the ideal. So it is quite nontrivial to prove that $\M$ is {\it globally\1} isomorphic to the nearby cycle $\D_X$-module $(V^0\!/V^1)\B_f(*Y)$. Note that Steenbrink's calculation only shows that its de~Rham complex ${\rm DR}_X(\M)$ is isomorphic to the nearby cycle complex $\psi_f\C_X[n]$ {\it after\1} taking the {\it cohomology sheaves.} Here it is not difficult to prove the isomorphisms between the {\it graded pieces\1} of the monodromy filtration. To compare the extension classes between them for the comparison of the $\ddd_1$-{\it differentials\1} of the two {\it weight spectral sequences,} one can apply the adjunction relations for the inclusions of strata to conclude that the extension classes are determined {\it locally.} This implies that the extension classes are given by the \v Cech-type restriction and Gysin morphisms with complex coefficients.
\sk
In the {\it reduced\1} case, the situation can be improved by considering the corresponding {\it left\1} $\D_X$-module, since we have the canonical generator $f^{-1}\!f^s\ins\OO_X(*Y)[s]f^s$, which is annihilated by the corresponding {\it left\1} ideal, where $s$ is identified with $-\dd_tt$. There is no canonical generator in the {\it non-reduced\1} case because of the difference between ``lower" and ``upper" extensions (corresponding to the intervals $(-1,0]$ and $[0,1)$ to which the eigenvalues of the residue belong) unless $Y_{\rm red}\sst X$ is defined by a global function on $X$ (shrinking $S$ if necessary). We can however apply Remark~\hl{R2.2}{2.2} for the proof of the isomorphisms
\htt{2.3.4}{}
$$\aligned{\rm DR}_X^{-1}\bl(\Om_{X/S}^{\ssb}(\log Y)(-Y_{\rm red})/\Om_{X/S}^{\ssb}(\log Y)(-Y)\br)&=\D_Xf^s/V^1\D_Xf^s,\\{\rm DR}_X^{-1}\bl(\Om_{X/S}^{\ssb}(\log Y)/\Om_{X/S}^{\ssb}(\log Y)(-Y_{\rm red})\br)&=V^1\D_Xf^s/\D_Xf^{s+1},\endaligned
\leqno(2.3.4)$$
where $\D_Xf^s$ is identified with $V^{>0}\B_f$ (and $f^s$ with $\delta(f{-}t)$).
\par\htt{R2.3b}{}\msn
{\bf Remark~2.3b.} It is well known that the {\it local invariant cycle theorem\1} is equivalent to the {\it decomposition theorem\1} over a curve, see for instance Appendix in \cite{KLS} and \cite{pk} for the decomposition theorem of the direct images of constant sheaf by proper morphisms from K\"ahler manifolds. (These are related to \cite{Ch}.) Note also that the local invariant cycle theorem holds in the case where the total space $X$ is nonsingular and K\"ahler and $f$ is proper {\it without\1} assuming that the singular fiber is a {\it divisor with normal crossings.} Indeed, the constant sheaf $\Q_X$ is a direct factor of $\RR\pi_*\Q_{\Xt}$ with $\pi\col\Xt\tos X$ an embedded resolution of $Y$ using only {\it Verdier duality.} The K\"ahler condition in the normal crossing case can be replaced by the condition written before Remark~\hl{R2.3a}{2.3a}.
\par\htt{R2.3c}{}\msn
{\bf Remark~2.3c.} The sign system in \cite[2.3]{Ch} adopts Deligne's one in \cite{GN}, where ``positivity" of polarization of Hodge structure is different from the usual one. Indeed, $(-1)^qS(-,-)$ is used there, and this is compatible with the sign in \cite[(2.1.15)]{De2}, since $(2\pi i)^wi^{q-p}\eq (-1)^q(2\pi)^w$ if $p\pl q\eq w$.
\par\htt{R2.3d}{}\msn
{\bf Remark~2.3d.} Concerning the {\it complex\1} Mellin transformation, we have
\htt{2.3.5}{}
$$\int_{\De_r}t^j\overline{t}{}^k|t|^{2s+2\al}\,\tfrac{\ddd t}{t}{\wedge}\tfrac{\ddd\overline{t}}{\overline{t}}\eq{-}2\pi i\,\delta_{j,k}\,r^{2s+2\al+2j}\1(s{+}\al{+}j)^{-1}\q(j,k\ins\N),
\leqno(2.3.5)$$
and its residue at $s\eq{-}\al$ is equal to $2\pi\delta_{j,0}\delta_{k,0}/i$, which is independent of $r\sgt0$, where $\De_r$ is a disk of radius $r$ with coordinate $t$. This is closely related to the formulas before \cite[Ex.\,6.1]{Ch} and at the end of the proof of \cite[Prop.\,5.1.14]{Sab2}.
\par\htt{R2.3e}{}\msn
{\bf Remark~2.3e.} Related to the {\it trace morphism,} let $E=\C/(\Z\pl\Z\tau)$ with ${\rm Im}\,\tau\sgt0$ choosing $i\eq\sqrt{-1}$. Set $\xi\defs\ddd z\ins\Gamma(E,\Om_E^1)$ with $z$ the coordinate of $\C$. Let $u,u'$ be the dual basis of the basis of $H_1(E,\Z)$ consisting of the images of the segments $I$ and $I'$ between $0$ and $1$ and $0$ and $\tau$ respectively. Since $\int_I\ddd z\eq1$, $\int_{I'}\ddd z\eq\tau$, we get that $\xi\eq u\pl\tau u'$, hence $\xib\eq u\pl\taub u'$, and
\htt{2.3.6}{}
$$\xi{\wedge}\1\xib\eq(\taub{-}\tau)u{\wedge}u'\eq{-}2i\,{\rm Im}\,\tau\,u{\wedge}u'.
\leqno(2.3.6)$$
Note that $u{\wedge}u'\ins H^2(E,\Z)$ is the generator corresponding to the orientation given by the choice $i\eq\sqrt{-1}$. On the other hand, setting $z\eq x\pl iy$, we have
\htt{2.3.7}{}
$$\int_E\xi{\wedge}\1\xib\eq\int_{\Xi}\ddd z{\wedge}\ddd\zb\eq{-}2i\!\int_{\Xi}\ddd x{\wedge}\ddd y\eq{-}2i\,{\rm Im}\,\tau,
\leqno(2.3.7)$$
where $\Xi\sst\C$ is the parallelogram with vertices $0,1,\tau,\tau{+}1$.
\sk
These show that the analytic trace morphism coincides with the topological one {\it without\1} multiplication by $2\pi i$ in this case. If we change the choice $i\eq\sqrt{-1}$, then $\tau$ must be changed by $-\tau$, but ${\rm Im}\,\tau$ remains {\it unchanged.}

\sk
{\smaller\smaller 
RIMS Kyoto University, Kyoto 606-8502 Japan}
\end{document}